
\input amstex
\documentstyle{amsppt}
\hsize = 5.4 truein
\vsize = 8.7 truein
\baselineskip = 12pt
\NoBlackBoxes
\TagsAsMath

\define\si{\sigma}

\define\orht{\overset\rightharpoonup\to}

\define\la{\lambda}
\define\k{\kappa}

\define\a{\alpha}
\define\be{\beta}
\define\de{\delta}
\define\De{\Delta}
\define\ga{\gamma}
\define\al{\aleph}

\define\th{\theta}

\define\om{\omega}

\define\z{\zeta}

\define\pr{\prime}

\define\1{\bigskip}
\topmatter
\title Filters, Cohen Sets and Consistent Extensions\\
of the Erd\"os-Dushnik-Miller Theorem
\endtitle
\author Saharon Shelah and Lee J. Stanley
\endauthor
\thanks {The research of the
first author was partially supported by the NSF and
the Basic Research Fund, Israel Academy of Science.  This is
paper number 419 in the first author's list of publications.
Theorem 1 was proved in Fall, 1989 when both authors benefitted
from the hospitality of MSRI, for which they record their gratitude.
Theorem 3 was proved in Fall, 1993; Theorem 4 was
was proved in Fall 1994.}
\endthanks
\address{Hebrew University, Rutgers University}
\endaddress
\address{Lehigh University}
\endaddress

\abstract{We present two different types of models
where, for certain singular cardinals $\la$ of uncountable
cofinality, $\la \rightarrow (\la, \om + 1)^2$, although $\la$
is not a strong limit cardinal.  We announce, here, and will
present in a subsequent paper, \cite{7}, that, for example,
consistently, $\al_{\om_1}
\not\rightarrow (\al_{\om_1},\ \om + 1)^2$ and
consistently, $2^{\al_0} \not\rightarrow (2^{\al_0}, \om + 1)^2$.}
\endabstract

\endtopmatter
\1\1
\subheading{\S 0.  INTRODUCTION}
\1
For regular uncountable $\k$, the Erd\"os-Dushnik-Miller theorem,
Theorem 11.3 of \cite{2}, states that
$\k \rightarrow (\k,\ \om + 1)^2$.
For singular cardinals, $\k$, they were only able to obtain
the weaker result, Theorem 11.1 of \cite{1}, that
$\k \rightarrow (\k,\ \om)^2$.
It is not hard to see that if $cf\ \k = \om$ then
$\k \not\rightarrow (\k,\ \om + 1)^2$.
If $cf\ \k > \om$ and
$\k$ is a strong limit cardinal, then it follows from the
General Canonization Lemma, Lemma 28.1 of \cite{1}, that
$\k \rightarrow (\k,\ \om + 1)^2$.  Question 11.4 of \cite{1}
is whether this holds without the assumption that $\k$ is a strong limit
cardinal, e.g., whether,
in ZFC,

$$ (1)\ \ \al_{\om_1} \rightarrow (\al_{\om_1},\ \om + 1)^2.$$

Another natural question, which the second author first heard from
Todorcevic, is whether, in ZFC,

$$ (2)\ \ 2^{\al_0} \rightarrow
(2^{\al_0},\ \om + 1)^2.$$

In connection with (2), we note that the first author proved,
\cite{2}, \S 2, the consistency of $2^{\al_0} \rightarrow
[\al_1]^2_{n,2}$.

In this paper we address these questions, by presenting
two types of models where there is a singular
cardinal $\la$ of uncountable cofinality, such that $\la \rightarrow
(\la,\ \om + 1)^2$ even though $\la$ {\it is not a strong limit cardinal}.
In either model, $\la$ can be taken to be $\al_{\om_1}$ and in the
second, we can also have, simultaneously, $\la = 2^{\al_0}$.
We also announce here, and will present in a subsequent paper, some
very recent results that show that, consistently,
(1) and (2) above may fail.  For (1), this answers
Question 11.4 of \cite{1} negatively.

The first type of model seems specific to
having the order type of the homogeneous
set for the second color (green, for us, whereas the first color is
the \lq\lq traditional" red)
be $\om + 1$,
whereas the second model allows
generalizations to green homogeneous sets of
order type $\th + 1$ for cardinals, $\th$, with $\om \leq \th < cf\ \la$,
under appropriate hypotheses.
On the other hand, the proof for the first model is an outright implication
from a hypothesis which follows from the existence of certain
partition cardinals, either outright, or in inner models,
and therefore, certainly,
from the failure of the SCH, for example.
\1\1
\proclaim{Theorem 1} If $\om < \k = cf\ \la,\ 2^\k < \la$ and
there is a normal {\it nice filter} on $\k$, then $\la \rightarrow
(\la,\ \om + 1)^2$.
\endproclaim
\1

There is no assumption on powersets between $\k$ and $\la$.  We prove
Theorem 1 in \S 1.  The notion of {\it nice filter}
is due to the first author.
In (1.1), below, we will give a condensed definition,
sufficient for our purposes, which is
consistent with the more general treatment of \S\S 0,1 of Chapter
V of \cite{3}.  This is essentially clause (2) of
Definition V.1.9 of \cite{3}.  The crucial property of nice filters,
for the purposes of this paper,
is that we can define a certain kind of rank function,
$rk^2_D(f,\ \Cal E)$, with ordinal values,
where $D$ is any normal nice filter,
$f:\k \rightarrow OR$ and $\Cal E$ is the family
of normal nice filters on $\k$.  This rank function has the following
important property:
\1\1
\roster
\item"{(\#)}" If $D \in \Cal E,\ f,\ g :\k \rightarrow OR,\ X$
is $D\text{-positive}$
and for all $\ga \in X,\ g(\ga) < f(\ga)$, then
then there is $D^\pr \in \Cal E$ with
$D \cup \{X\} \subseteq D^\pr$ such that $rk^2_{D^\pr}(g,\ \Cal E) <
rk^2_D(f,\ \Cal E)$.
\endroster
\1\1
\noindent
This can be extracted from the following items of  Chapter V of \cite{3}:
Claim V.2.13, and clause (1) of Fact V.3.16.
The existence of a nice filter on
$\om_1$, for example, is an outright consequence of the existence of a
$\mu$ such that $\mu \rightarrow (\a)^{<\om}_{\al_o}$ for all
$\a < (2^{2^{\al_1}})^+$.  It can also be obtained in forcing
extensions starting from models with such large cardinals.
For these results, see
Conclusion V.1.13 and Remark V.1.13A of \cite{3}.  In view
of the first fact, we easily have the following corollary to Theorem 1;
a later result in a similar vein is Woodin's striking result that
from $CH$ and the existence of a measurable cardinal it follows
that the club filter on $\al_1$ is not $\al_2\text{-saturated.}$
\1\1
\noindent
\proclaim{Corollary 2}  Assume that there is a measurable cardinal and
that $2^{\al_1} < \al_{\om_1}$.
Then $\aleph_{\omega_1} \rightarrow (\al_{\om_1}, \om + 1)^2$.
\endproclaim
\1\1

In the second type of model, we have several parameters.
We let $\om < \k = cf\ \la < \la$.
As mentioned above, we have a cardinal $\th$ with
$\om \leq \th < \k$.  We have an additional
cardinal parameter, $\si$, with $\si \neq \k$ and
$\th \leq \si$.  The cases $\si < \k$ and  $\si > \k$
require somewhat different different treatment, and lead
to Theorems 3 and 4, below, respectively, proved in \S 2 and
\S 3.  However, much of the preliminary material developed
for Theorem 3 carries over to the proof of Theorem 4.
The main case of Theorem 3 is when $\th = \si$, and the
connection to Theorem 1 is when $\th = \si = \omega$.
Theorem 3 was proved in Fall 1993 and Theorem 4 was proved
in Fall  1994.

For both Theorems, we assume that in $V,\ \la$ is a strong limit cardinal,
and that $\si^{<\si} = \si$.  Our model is obtained by
forcing with $\bold P$, which is the partial ordering for adding
at least $\la$ Cohen subsets of
$\si$.
When $\th > \om$, we need additional assumptions
to guarantee, for example, that in $V^{\bold P},\
\k \rightarrow (\k,\ \th + 1)^2$.  When $\th = \om$,
this is just the Erd\"os-Dushnik-Miller theorem for
$\k$.  The additional assumptions will involve
cardinal exponentiation, and will be discussed
below.  We then have:
\1\1
\proclaim{Theorem 3}  Suppose that in $V,\ \om \leq
\th \leq \si = \si^{<\si} <
\k = cf\ \la < \la \leq \nu,\ \la$ is a strong limit cardinal and for
all $\mu < \k,\ \mu^{<\th} < \k$.  Let $\bold P$ be the partial
ordering for adding $\nu$ Cohen subsets of $\si$.  Then, in
$V^{\bold P},\ \la \rightarrow (\la,\ \th + 1)^2$.
\endproclaim
\1\1
\proclaim{Theorem 4}  Suppose that in $V,\ \om \leq
\th < \k <  \si = \si^{<\si} < \la,\ \k = cf\ \la < \la \leq \nu,\
\la$ is a strong limit cardinal and for
all $\mu < \k,\ \mu^{<\th} < \k$.  Let $\bold P$ be the partial
ordering for adding $\nu$ Cohen subsets of $\si$.  Then, in
$V^{\bold P},\ \la \rightarrow (\la,\ \th + 1)^2$.
\endproclaim
\1\1
We shall deduce Theorem 4 from the following result about lifting
certain positive partition relations on $\k$
in $V$ to $\la$ in models,  $V^{\bold P}$, where $\k,\ \la,\ \si,\
\bold P$ are as in Theorem 4.
\1\1
\proclaim{Theorem 4*}  Suppose that in  $V,\ \om < \k < \si =
\si^{<\si} < \la,\ \k = cf\ \la \leq \nu,\ \la$ is a strong limit cardinal
and $\bold P$ is the partial ordering for adding $\nu$ Cohen subsets of
$\si$.  Suppose, further, that $\zeta < \k$ and that $\k \rightarrow
(\k,\ \zeta)^2$.  Then, in $V^{\bold P},\ \la \rightarrow
(\la,\ \zeta)^2$.
\endproclaim
\1\1
Of course, when we invoke Theorem 4* to obtain Theorem 4, we shall
take  $\zeta = \theta + 1$, and we will use the additional
hypotheses on cardinal exponentiation in $V$ to obtain the
hypothesis of Theorem 4*, that $\k \rightarrow (\k,\ \zeta)^2$.
Then, this relation will also hold in $V^{\bold P}$, since there
are no  new subsets of $\k$.  In fact, it is even possible to
factor Theorem 3 through a similar kind of result about lifting
positive relations on $\k$ to $\la$, but now lifting a  $V^{\bold P}$
relation on  $\k$ to a $V^{\bold P}$ relation on $\la$, since
this time, forcing with $\bold P$  will not necessarily preserve
a positive $V$ relation on $\k$.
In what follows, we shall not proceed in this fashion; however, we
do state the lifting theorem:
\1\1
\proclaim{Theorem 3*}  Suppose that in  $V,\ \zeta \leq \si + 1,
\ \si = \si^{<\si} < \k = cf\ \la \leq \nu,\ \la$ is a strong limit cardinal
and $\bold P$ is the partial ordering for adding $\nu$ Cohen subsets of
$\si$.  Suppose, further, that in $V^{\bold P},\ \k \rightarrow
(\k,\ \zeta)^2$.  Then, in $V^{\bold P},\ \la \rightarrow
(\la,\ \zeta)^2$.
\endproclaim
\1\1

Once again, in order to obtain Theorem 3 from
Theorem 3*, the additional
hypotheses in Theorem 3
on cardinal exponentiation in  $V$ are designed to
guarantee that the needed positive relation does hold
in $V^{\bold P}$.  It would, of course, be possible to combine
Theorems 3* and 4* into a single statement, but the proof would
certainly reflect the division into cases, which, here, is transparent
in the statements.
\1\1
Finally, though these more recent results
will be presented in a subsequent paper, \cite{7}, we state
here, as numbered theorems, the negative consistency
results for questions (1) and (2), mentioned above and
in the Abstract.
\1\1
\proclaim{Theorem 5}  Suppose that, in $L,\ \mu \geq \la > \k$ are cardinals,
$cf(\la) = \k > \om$
(for example, $\la = (\al_\k)^L,\ \k = (\al_1)^L$).  Let $G$ be
$\bold P\text{-generic over } L$, where $\bold P$ is
($L's$ version of) the partial
order for adding $\mu$ Cohen subsets of $\k$.  Then, in $L[G],\
\la \rightarrow (\la,\ \om + 1)^2$ iff, in $L,\ \k$ is weakly compact.
\endproclaim
\1\1
Taking $\k = (\al_1)^L,\ \la = (\al_\k)^L$, we get the negative consistency
result for (1).  Combining the methods used to obtain Theorem 5 for
this choice of $\k$ and $\la$, an additional forcing to add $\la$ Cohen
reals, and a double $\Delta\text{-system}$ argument for the second
forcing, we get:
\1\1
\proclaim{Theorem 6}  Con(ZFC) implies Con(ZFC \& $2^{\al_0}
\not\rightarrow (2^{\al_0}, \om + 1)^2$).
\endproclaim
\1\1
\proclaim{Remarks}
\endproclaim
\medskip
\roster
\item The proof of the Erd\"os-Dushnik-Miller theorem proceeds by assuming
that there is no homogeneous set of power $\k$ for the first color
(red, for us), and showing that a certain tree of homogeneous green sets
must have a branch of length $\om + 1$, which naturally yields a homogeneous
green set of order type $\om + 1$.  If $\th > \om,\ \th$ is a cardinal,
$\tau > \th$ is regular, and if:

$$(\ast)\text{  for all } \nu < \tau,\ \nu^{<\th} < \tau,$$

\noindent
then we can carry out essentially the same proof to show that
$\tau \rightarrow (\tau,\ \th + 1)^2$.
Thus, taking $\tau = \k$,
our hypotheses on cardinal exponentiation
in Theorem 3, which remain true in $V^{\bold P}$, do guarantee that
in $V^{\bold P},\ \k \rightarrow (\k,\ \th + 1)^2$.

Similary, if $\om < \zeta,\ \th = card\ \zeta,\ \th < \tau, \ \tau$
is regular and if:

$$(\ast\ast)\text{ for all } \nu < \tau, \ \nu^\th < \tau,$$

\noindent
then a similar tree argument shows that $\tau \rightarrow (\tau,\ \zeta)^2$.
Thus, the additional hypotheses on cardinal exponentiation in Theorem 4
do guarantee that we have the hypotheses of Theorem 4*.

For both theorems, we will also
need to know that for many successor cardinals, $\tau$, between
$\k$ and $\la$, we will have $\tau \rightarrow (\tau,\ \th + 1)^2$,
or $\tau \rightarrow (\tau, \zeta)^2$ (for Theorem 4*).
In view of the preceding paragraphs, it will suffice to have $(\ast)$
or $(\ast\ast)$ for $\tau$, in $V$.

One way of achieving this is to appeal to the fact that, in $V$,
$\la$ is a strong limit cardinal, and, for example, to take
$\tau = \mu^+$, where $\mu = \mu^\th$, and where $\mu$ is chosen
to have various other properties, as desired.
\medskip
\item  In all of what follows we shall have $\om < \k = cf\ \la < \la$.
We shall express $\la$ as  $sup\{\la_\eta|\eta < \k\}$, where
$(\la_\eta|\eta < \k)$ is increasing and continuous, and for $\eta = 0$
or $\eta$ a successor ordinal, $\la_\eta$ is a successor cardinal.
Various other properties of the $\la_\eta$ for such $\eta$ will
be introduced as needed.  One such property will be that
$\la_\eta = \mu^+$, where $\mu = \mu^\th$ (and has various other
properties, as desired).
We also let $\De_0 = \la_0$ and
for $\eta < \k,\ \De_{1 + \eta} = [\la_\eta,\ \la_{\eta + 1})$.  For
$\a < \la$ we will let $\eta(\a) = $ the unique $\eta < \k$ such
that $\a \in \De_\eta$.
\item  Investigation of the case $\si = \k$, which is not
treated in this paper, led to Theorems 5, 6, above, among
other results.  When $\si < \k$,
we use the $\si^+\text{-chain condition of } \bold P$, whereas
when $\si > \k$, we use the $<\si\text{-completeness}$ of $\bold P$.
\item  In Theorems 3* and 4*, it is clearly necessary that in
$V^{\bold P},\ \k \rightarrow (\k,\ \th + 1)^2$, respectively,
that $\k \rightarrow (\k,\ \zeta)^2$.
\item Our notation and terminology is intended
to either be standard or have a clear meaning, e.g., $card\ X$ for
the cardinality of $X$, $o.t.\ X$ for the order type of $X$, etc.
\item  Theorems 3 and 5 of \cite{6} are close in spirit to
some of the above material.  There are also similarities to certain
themes from \cite{4} and \cite{5}.
\endroster
\subheading{\S 1.  USING NICE FILTERS}
\1
In this section we prove Theorem 1 of the Introduction, which, for
convenience, we now restate.
\1\1
\proclaim{Theorem 1} If $\om < \k = cf\ \la,\ 2^\k < \la$ and
there is a normal {\it nice filter} on $\k$, then $\la \rightarrow
(\la,\ \om + 1)^2$.
\endproclaim
\1\1
\demo{Proof}
We begin by providing the promised definition of nice filter on $\k$.
If $D$ is a normal filter on $\k$ and $g$ is an ordinal valued
function with domain $\k$, we first define the game $Gw^*(D,\ g)$, as
follows.  On move $0$, player I chooses $D_0\ :=\ D$, and player II chooses
$A_0 \in (D_0)^+$, and chooses $g_0\ :=\ g$.  On move $n + 1$,
player I chooses $D_{n + 1}$, a normal filter on $\k$ extending
$D_n \cup \{ A_n \}$, and player II chooses $A_{n + 1} \in (D_{n + 1})^+$,
AND $g_{n + 1} <_{D_{n + 1}^*} g_n$, where $D_{n + 1}^*\ :=\ $
the normal filter on $\k$ generated by $D_{n + 1} \cup \{ A_{n + 1} \}$.
Player I wins if at some stage $n + 1$, Player II has no legal play.
We then state:
\1\1
\proclaim{(1.1)  Definition}  $D$ is nice if for all ordinal valued
functions, $g$, with domain $\k$, Player I has a winning strategy
in $Gw^*(D,\ g)$.
\endproclaim
\1\1
Proceeding with the proof of the Theorem, we assume
that $\om < \k = cf\ \la < \la,\ 2^\k < \la$ and that there is a
nice normal filter on $\k$.
We will show that $\la \rightarrow (\la,\ \om + 1)^2$.
There are no assumptions about powers of cardinals larger than $\k$, and,
as noted in the Introduction, the interest of the result is when $\la$
is not a strong limit cardinal.  The simplest case, of course, is when
$\k = \al_1$ and $\la =\al_\k$.

So, towards a contradiction, suppose that $c:[\la]^2\ \rightarrow\
\{\text{red, green}\}$ but has no red set of power $\la$ and no
green set of order type $\om + 1$.  Let $\la_\eta,\ \De_\eta,\ \eta < \k$
be as in Remark 2 of the Introduction.  We can clearly assume, in
addition, that $\la_0 > 2^\k$, for
$\eta < \k,\ \la_{\eta + 1} \geq \la_\eta^{++}$, and
that each $\De_\eta$ is homogeneous red for $c$.
The last is by the Erd\"os-Dushnik-Miller theorem for
$\la_{\eta + 1}$.

For $0 < \eta < \k$, we define $Seq_\eta$ to be
$\{(i_0, ..., i_{n-1})|\eta(i_0) < ... < \eta(i_{n-1}) < \eta\}$.  For $\z
\in \De_\eta$ and
$(i_0, ..., i_{n-1}) = \orht{i} \in Seq_\eta$, we say $\orht{i}
\in T^\z$ iff $\{i_0, ..., i_{n-1},\ \z\}$ is homogeneous green for $c$.
Note that an infinite decreasing (for reverse inclusion) branch
in $T^\z$ violates the nonexistence of a green set of order type
$\om + 1$, so, under reverse inclusion, $T^\z$ is well-founded.
Therefore the following definition of a rank function,
$rk^\z$, on $Seq_\eta$
can be carried out.

We define $rk^\z:Seq_\eta \rightarrow OR \cup \{-1\}$ by setting
$rk^\z(\orht{i})$ to be $-1$ if
$\orht{i}^{\frown}\z$ is not homogeneous green; otherwise, define
$rk^\z(\orht{i}) \geq \eta$ iff for all $\tau < \eta$ there is $j$
such that $rk^\z(\orht{i}^{\frown}j) \geq \tau$.  Of course, for limit
ordinals, $\de$, if for all $\eta < \de,\ rk^\z(\orht{i}) \geq \eta$,
then $rk^\z(\orht{i}) \geq \de$, and so for all $\orht{i} \in T^\z$,
there is a largest $\eta$ such that $rk^\z(\orht{i}) \geq \eta$.
We take $rk^\z(\orht{i})$ to be this largest $\eta$.  In fact, it is clear
that the range of $rk^\z$ is a proper initial segment of $\mu_\eta^+$, where
$\mu_\eta = card\ \bigcup\{\De_\tau|\tau < \eta\}$, and so, in particular,
the range of $rk^\z$ has power at most $\la_\eta$.  Note that $\la_{\eta + 1}
> \mu_\eta^+$.

But then,
we can find $B_\eta$ an end-segment of $\De_\eta$ such that for all
$\orht{i} \in Seq_\eta$ and all $0 \leq \ga < \mu_\eta^+$,
if there is $\z \in B_\eta$ such that $rk^\z(\orht{i}) = \ga$, then
there are $\la_{\eta + 1}$ such $\z$.  Recall that $\De_\eta$ and therefore
also $B_\eta$ are of order type $\la_{\eta + 1}$, which is a successor
cardinal.  Everything is now in place for the main definition.

\1\1

\proclaim{(1.2) Definition}  $(\orht{i},\ Z,\ D,\ f) \in K$ iff
\1
\roster
\item $D$ is a {\it nice}, normal filter on $\k$,
\item $f:\k \rightarrow OR$,
\item $Z \in D$,
\item for some $0 < \eta < \k,\ \orht{i} \in Seq_\eta$, and for all
$\tau \in Z \setminus (\eta + 1)$, there is $\z \in B_\tau$ such that
$rk^\z(\orht{i}) = f(\tau)$ (so, in particular, $\orht{i} \in T^\z$).
\endroster
\endproclaim

\1\1

Note that $K \neq \emptyset$, since if we choose $\z_\tau \in
B_\tau$, for $\tau < \k$, take $Z = \k,\ \orht{i} = $
the empty sequence, choose $D$ to be any nice normal filter on
$\k$ and define $f$ by $f(\tau) = rk^{\z_\tau}(\orht{i})$,
then $(\orht{i},\ Z,\ D,\ f) \in K$.

Now, let $\Cal E$ be the family of nice normal filters
on $\k$.  Since $rk^2_D(f,\ \Cal E) \in OR$, clearly
among the $(\orht{i},\ Z,\ D,\ f) \in K,
\text{ there is one with } rk^2_D(f,\ \Cal E) \text{ minimal.}$

So, fix one such, and denote it by $(\orht{i}^*,\ Z^*,\ D^*,\ f^*)$.
For $\tau \in Z^*$, set $C_\tau = \{\z \in B_\tau|rk^\z(\orht{i}) \leq
f^*(\tau)\}$.  Thus $card\ C_\tau = \la_{\tau + 1}$, and
for all $\z \in C_\tau,\ range(\orht{i}^* \cup \{\z\})$ is homogeneous green.
Now suppose $\tau \in Z^*$.  For all $\ga \in Z^* \setminus (\tau + 1)$
and $\z \in C_\tau$,
let $C^+_\ga (\z) = \{\xi \in C_\ga |c(\{\z,\xi\}) = \text{green}\}$.
Also, let $Z^+(\z) = \{\ga \in Z^* \setminus (\be + 1)| card\ C^+_\ga (\z)
= \la_{\ga + 1}\}$.  It is, perhaps, worth pointing out that we
could just as well
have required only that $C^+_\ga \neq \emptyset$.

\1

\proclaim{(1.3)  Lemma}  For a $D\text{-positive}$ set of
$\tau \in Z^*$ and for $\la_{\tau + 1}$ many $\z \in C_\tau,\
Z^+(\z)$ is $D\text{-positive}$.
\endproclaim
\medskip
\demo{Proof}  For $\tau \in Z^*$ and $\z \in C_\tau$, let
$Y(\z) = \k \setminus Z^+(\z)$.  Since $\la_0 > 2^\k$,
for all $\tau \in Z^*$ there is
$Y = Y_\tau \subseteq \k$ and $C^\pr_\tau \subseteq C_\tau$ with
$card\ C^\pr_\be = \la_{\tau + 1}$ such that for all $\z \in C^\pr_\tau,\
Y(\z) = Y_\be$.

Let $\hat Z = \{\tau \in Z^*|Y_\tau \in D\}$.
We now conclude by showing that $\hat Z \not\in D$.
If $\hat Z \in D$, then, since $D$ is normal,
we would have $Y^* \in D$, where
$Y^* = \{\tau \in \hat Z|\text{for all } \eta \in \hat Z \cap \tau,\
\tau \in Y_\eta\}$.
But then, by shrinking the $C^\pr_\tau$ for $\tau \in Y^*$, as
in the next paragraph,
we would get a homogeneous red set of power $\la$, which is impossible.

We define $\hat C_\tau$ for $\tau \in Y^*$ by recursion on $\tau$ in such
a way that
$\hat C_\tau$ is a subset of $C^\pr_\tau$ of power $\la_{\tau + 1}$.
So, let $\tau \in Y^*$, and set $\xi \in \hat C_\tau \text{ iff }
\xi \in C^\pr_\tau$
and for all $\eta \in Y^* \cap \tau$ and all $\z \in \hat C_\eta ,\ \xi
\not\in C^+_\tau (\z)$.  So, in fact, $\hat C_\tau$ is the result of
removing at most $\la_\tau$ elements from $C^\pr_\tau$.  But then, clearly
the union of the $\hat C_\tau$ for $\tau \in Y^*$ is homogeneous red.
This concludes the proof of Lemma 1.2.
\enddemo

\1

We maintain the notation of the proof of Lemma 1.2.
Fix $\tau$ as guaranteed by Lemma 1.2, i.e., such that
$Y_\tau$ is defined, but $Y_\tau \not\in D$.
Let $X = Z^* \setminus Y_\tau$.  Note that, for any $\z \in C^\pr_\tau,\
X \setminus (\tau + 1) = Z^+(\z)$ and $X$ is $D\text{-positive}$.
Now fix $\z \in C^\pr_\tau$.  For $\ga \in X \setminus (\tau + 1)$,
note that by the definition of $C^+_\ga (\z)$, there is $j \in C^+_\ga(\z)$
such that $rk^j(\orht{i}^*) \leq f^*(\ga)$.  Choose one such and call it
$j_\ga$.   Thus, again by the definition of $C^+_\ga(\z),\
\orht{i}^{*\frown}\z^{\frown}j_\ga$ is homogeneous green, and so, by the
definition of $rk^{j_\ga},\
rk^{j_\ga}(\orht{i}^{*\frown}\z) < f^*(\ga)$.

Now, define $g:\k \rightarrow OR$ by $g(\ga) =
rk^{j_\ga}(\orht{i}^{*\frown}\z)$, if $\ga \in X \setminus (\tau + 1)$,
and $g(\ga) = 0$, otherwise.  Now, by the definition of
$rk^2_D(f^*,\ \Cal E)$, (again, see Chapter 5 of \cite{3}) there is
$D^\pr \in \Cal E$ with $D \cup \{X\} \subseteq D^\pr$ and such that
$rk^2_{D^\pr}(g,\ \Cal E ) < rk^2_D (f^*,\ \Cal E )$.  However, it is
easily verified that $(\orht{i}^{*\frown}\z,\ X,\ D^\pr,\ g) \in K$, and,
finally, this contradicts the choice of
\newline
$(\orht{i}^*,\ Z^*,\ D,\ f^*)$,
and thus completes the proof of Theorem 1.
\enddemo
\1\1
\subheading{\S 2.  ADDING COHEN SETS BELOW THE COFINALITY}
\1
In this section, we prove Theorem 3 of the Introduction, whose
statement we now recall for convenience.
\1\1
\proclaim{Theorem 3}  Suppose that in $V,\ \om \leq
\th \leq \si = \si^{<\si} <
\k = cf\ \la < \la \leq \nu,\ \la$ is a strong limit cardinal and for
all $\mu < \k,\ \mu^{<\th} < \k$.  Let $\bold P$ be the partial
ordering for adding $\nu$ Cohen subsets of $\si$.  Then, in
$V^{\bold P},\ \la \rightarrow (\la,\ \th + 1)^2$.
\endproclaim
\1
\demo{Proof}
So, let $\la,\ \k,\ \th,\ \si,\ \nu,\ \bold P$
be as in the statement
of Theorem 3, and let $(\la_\eta|\eta < \k)$ be as in
Remark (2) of the Introduction, and suppose, in addition that
$\la_0 = \mu_0^+$, where $\mu_0 = (\beth_\om(\k))^\th$,
and for $\eta < \k,\ \la_{\eta + 1} = \mu_\eta^+$, where
$\mu_\eta = (\mu_\eta)^{((\beth_\om(\la_\eta))^\th)}$.  Thus, by Remark 1
of the Introduction, we will have that in $V^{\bold P},\
\k \rightarrow (\k,\ \th + 1)^2$ and similarly:

$$\align (!)\ \  &\text{ in } V^{\bold P}, \text{ for each }
\eta < \k \text{ which is either } 0\\
&\text{ or a successor ordinal, }
\la_\eta \rightarrow (\la_\eta,\ \th + 1)^2.\endalign$$

This follows from our choice of the $\la_\eta$ since
forcing with $\bold P$ adds no new sequences of ordinals of
length $< \th$.
Also, let $\De_\eta$, and $\eta(\a)$ be as in Remark 2 of
the Introduction.

For $A \subseteq \nu$, we let $\bold P|A$ be the subordering of
$\bold P$ with underlying set the set of $p \in P$ with domain
included in $A$.  If $card\ A = card\ B$ and $T$ is
a bijection from $A$ to $B$, we abuse
notation by also taking $T$ to be the isomorphism from $\bold P|A$
to $\bold P|B$ induced by $T$.

Suppose, now, that $\bold c$ is a
$\bold P\text{-name}$ and that $p \in P$ forces that
$\bold c :[\la]^2 \rightarrow
\{\text{red, green}\}$.  We now embark on an analysis of $\bold c$
as a $\bold P\text{-name}$ culminating in (*), following (2.9).
This analysis carries over to \S 3, and even in the case $\si = \k$.
We use the latter case in our forthcoming paper, \cite{7},
when $\k$ is weakly compact.  Therefore, we temporarily
drop the assumption the assumption $\si < \k$, or even that $\si \neq
\k$, retaining only that $\si = \si^{<\si} < \la$.

By (!), we can assume, without loss of generality,
that for each $\eta < \k,\ p$ forces that
$\De_\eta$ is
homogeneous red for $\bold c$.
In order to develop material that will carry over
to the proof of Theorem 4*, in \S 3, for now, we make
no additional hypotheses about $\bold c$.

For $\a < \be < \la$,
let $A(\a,\ \be)$ be a subset of $\nu$ of power at most $\si$ such
that $\bold c (\a,\ \be)$ is a $\bold P|A(\a,\ \be)\text{-name}$.
Such $A(\a,\ \be)$ exists, since $\bold P$ has the $\si^+\text{-cc}$.
Let $A^* = \bigcup\{A(\a,\ \be)|\a < \be < \la\}$ and let
$\bold P^* = \bold P|A^*$.  Without loss of generality,
$dom\ p \subseteq A^*$.  Thus, $\bold c \in V^{\bold P^*}$,
so by arguing in $V^{\bold P^*}$, and remarking that
$card\ A^* = \la$ and therefore that $\bold P^* \cong \bold P|\la$,
we can assume, without loss of generality, that $\nu = \la$, which
we do from here on.

For $\a < \be < \la$, let $\pi(\a,\ \be) = o.t.\ A(\a,\ \be)$ and
let $(\rho^{\a,\ \be}_\zeta|\zeta < \pi(\a,\ \be))$ be the
increasing enumeration of $A(\a,\ \be)$.  Also, let $T(\a,\ \be)$
be the order isomorphism from $A(\a,\ \be)$ to $\pi(\a,\ \be)$
(so $T(\a,\ \be)(\rho^{\a,\ \be}_\zeta) = \zeta$).  Let
$\bold{c^\pr}(\a,\ \be)$ be the $\bold P|\pi(\a,\ \be)\text{-name}$
which results from applying $T(\a,\ \be)$ to $\bold c(\a,\ \be)$ where
$T(\a,\ \be)$ is viewed as the isomorphism from $\bold P|A(\a,\ \be)$
to $\bold P|\pi(\a,\ \be)$, as in the previous paragraph.
Fix
functions $F_i:[\la]^2 \rightarrow \la$, for $i < \si$, such
that for $\a < \be < \la,\ A(\a,\ \be) = \{F_i(\a,\ \be)|i < \si\}$.

\proclaim{(2.1) Definition}
Let $Y(\a,\ \be) = \{(i,\ \zeta) \in \si\times\pi(\a,\ \be)|
F_i(\a,\ \be) = \rho^{\a,\ \be}_\zeta\}$.
We also
let $X$ be the set of ordered 4-tuples, $(\a,\ \be,\ \ga,\ \de)$
from $\la$ such that $\a < \be$ and $\ga < \de$, and
we define a function
$c^*$ with domain $X$ by:
\1
$c^*(\a,\ \be,\ \ga,\ \de) = (\pi(\a,\ \be),\ \pi(\ga,\ \de),\
Y(\a,\ \be),\ Y(\ga,\ \de),\ \bold{c^\pr}(\a,\ \be),\ \bold{c^\pr}(\ga,\ \de))$.
\1

Note that the following set is easily recoverable from
$c^*(\a,\ \be,\ \ga,\ \de)$:
\1
$$\hat c (\a,\ \be,\ \ga,\ \de) =
\{(i,\ j)|F_i(\a,\ \be) = F_j(\ga,\ \de)\}.$$

\noindent
We abuse notation below
by acting as if this were actually part of $c^*(\a,\ \be,\ \ga,\ \de)$.
Also note that $range\ c^*$ has power at most $2^\si$.
\endproclaim

Applying the general canonization lemma, Lemma 28.1 of \cite{1} to
$c^*$, we get $B_\eta \subseteq \De_\eta$ with $card\ B_0 > \k + \si$
and for $0 < \eta < \k,\ card\ B_\eta > \la_\eta$,
and such that $(B_\eta:\eta < \k)$ is canonical for $c^*$, i.e,
letting $B = \bigcup\{B_\eta|\eta < \k\}$, if
$(\a_n|n < 4),\ (\be_n|n < 4) \in X \cap B^4$
and for all $n < 4,\ \eta(\a_n) =
\eta(\be_n)$, then $c^*((\a_n|n < 4)) = c^*((\be_n|n < 4))$.

Further note that if $\eta(\a_1) = \eta(\a_2) < \eta(\be_1) =
\eta(\be_2)$ and $\a_1,\ \a_2,\ \be_1,\ \be_2 \in B$, then
since $c^*(\a_1,\ \be_1,\ \a_1,\ \be_1)
= c^*(\a_2,\ \be_2,\ \a_2,\ \be_2)$, we also have that
$\bold c^\pr(\a_1,\ \be_1) =
\bold c^\pr(\a_2,\ \be_2)$.  This, in turn,
means that if $p_1 \in P|A(\a_1,\ \be_1),\
p_2 = (T(\a_2,\ \be_2))^{-1}\circ T(\a_1,\ \be_1)(p_1)$,
and $x \in \{\text{red, green}\}$, then $p_1$ forces
that $\bold c (\a_1,\ \be_1) = x$ iff $p_2$ forces
that $\bold c (\a_2,\ \be_2) = x$.  We will use this
observation in several places in what follows.
\1
\proclaim{(2.2) Lemma}  Suppose that $(\a,\ \be,\ \ga,\ \de) \in X \cap
B^4,\ \a \not\in \{\ga,\ \de\},\ \a^\pr \in B_{\eta(\a)}$ and
$F_i(\a,\ \be) = F_j(\ga,\ \de)$.  Then also $F_i(\a^\pr,\ \be) =
F_j(\ga,\ \de)$, and analogous statements hold where the
values of the other coordinates of $(\a,\ \be,\ \ga,\ \de)$
are varied instead of varying the first coordinate.
\endproclaim

\demo{Proof}  This is clear since $(B_\eta|\eta < \k)$ is canonical
for $c^*,\ \eta(\a) = \eta(\a^\pr)$ and so, as noted
at the end of Definition 2.1, $(i,\ j) \in
\hat c (\a,\ \be,\ \ga,\ \de)$ iff $(i,\ j) \in
\hat c (\a^\pr,\ \be,\ \ga,\ \de)$.
\enddemo

\proclaim{(2.3) Definition}  Suppose $\eta < \tau < \k,\
i < \si,\ \a \in B_\eta,\ \be \in B_\tau$.
We define
$F_{i,\ \a}^\tau:B_\tau \rightarrow \la$ and
$F_{i,\ \eta}^\be:B_\eta \rightarrow \la$ by $F_{i,\ \a}^\tau(\be^\pr) =
F_i(\a,\ \be^\pr)$ and $F_{i,\ \eta}^\be(\a^\pr) = F_i(\a^\pr,\ \be)$.
\endproclaim
\1
\proclaim{(2.4) Lemma}  If $\eta < \tau < \k,\ i < \si$ then:

\1
\roster
\item either (for all $\a \in B_\eta,\ F_{i,\ \a}^\tau$
is constant) or (for all $\a \in B_\eta,\ F_{i,\ \a}^\tau$ is one-to-one).
\medskip
\item either (for all $\be \in B_\tau,\ F_{i,\ \eta}^\be$ is constant)
or (for all $\be \in B_\tau,\ F_{i,\ \eta}^\be$ is one-to-one).
\endroster
\endproclaim

\demo{Proof}  We first argue that each is either constant or one-to-one.
We consider the $F_{i,\ \a}^\tau$.  Let $\be_1 \neq \be_2$ both in
$B_\tau$.  We claim that if $F_i(\a,\ \be_1) =
F_i(\a,\ \be_2)$ then $F_{i,\ \a}^\tau$ is constant, while if
$F_i(\a,\ \be_1) \neq F_i(\a,\ \be_1)$, then $F_{i,\ \a}^\tau$ is
one-to-one.  In the first case, $(i,\ i) \in
\hat c (\a,\ \be_1,\ \a,\ \be_2)$,
while in the second case, $(i,\ i) \not\in
\hat c (\a,\ \be_1,\ \a,\ \be_2)$.
But then, by canonicity, if $\be \in B_\tau \setminus \{\be_1,\ \be_2\}$,
$(i,\ i) \in \hat c (\a,\ \be,\ \a,\ \be_1)$
iff $(i,\ i) \in \hat c (\a,\ \be,\ \a,\ \be_2)$ iff $(i,\ i) \in
\hat c (\a,\ \be_1,\ \a,\ \be_2)$.  If $(i,\ i)$ is a member of none, then
$F_{i,\ \a}$ is one-to-one.  If $(i,\ i)$ is a member of all, then
$F_{i,\ \a}$ is constant.
The argument for the $F_{i,\ \eta}^\be$ is completely analogous.

We now argue that if $\a_1 \neq \a_2$ both in $B_\eta$,
and $F_{i,\ \a_1}^\tau$ is constant then so is
$F_{i,\ \a_2}^\tau$.
Once again, the argument for $F_{i,\ \eta}^{\be_1}$ and
$F_{i,\ \eta}^{\be_2}$ is completely analogous.
So, suppose that $F_{i,\ \a_1}^\tau$ is constant.

Choose $\ga_1 \neq \ga_2$ both in $B_\tau$.  Since
$F_{i,\ \a_1}^\tau$ is constant, $(i,\ i) \in
\hat c (\a_1,\ \ga_1,\ \a_1,\ \ga_2)$, so, by canonicity,
$(i,\ i) \in \hat c (\a_2,\ \ga_1,\ \a_2,\ \ga_2)$ which
means that $F_{i,\ \a_2}^\tau$ is constant.
\enddemo
\1
\proclaim{(2.5) Remark}  In Lemma 2.4, we cannot conclude that
if $F_{i,\ \a}^\tau$ is constant (resp. one-to-one) then
$F_{i,\ \eta}^\be$ is constant (resp. one-to-one), as this
would involve an \lq\lq illegal" application of canonization,
comparing a \lq\lq 1,2" case to a \lq\lq 2,1" case.  It is,
however, worth noting that if all the $F_{i,\ \a}^\tau$ are constant,
then all the $F_{i,\ \eta}^\be$ are constant iff all the
$F_{i,\ \a}^\tau$ have the same constant value; similarly,
if all the $F_{i,\ \eta}^\be$ are constant, then all the
$F_{i,\ \a}^\tau$ are constant iff all the $F_{i,\ \eta}^\be$
have the same constant value.  We argue for the first statement.

Suppose that all the $F_{i,\ \a}^\tau$
are constant.  Let $\a_1 \neq \a_2$ both in $B_\eta$ and $\be \in B_\tau$.
Then $F_{i,\ \eta}^\be$ is constant iff $F_i(\a_1,\ \be) =
F_i(\a_2,\ \be)$ and therefore, since the $F_{i,\ \a_j}^\tau$ are constant,
this holds iff they have the same constant value.
\endproclaim
\1
\proclaim{(2.6) Definition}  For $\eta < \tau < \k,\ i < \si,\ \a \in B_\eta,\
\be \in B_\tau$, we define $F_i(\a,\ \tau),\ F_i(\eta,\ \be)$ by
$F_i(\a,\ \tau) = $ the constant value of $F_{i,\ \a}^\tau$,
if $F_{i,\ \a}^\tau$ is a constant function, and undefined if
it is a one-to-one function.  Similarly, $F_i(\eta,\ \be) = $ the constant
value of $F_{i,\ \eta}^\be$, if $F_{i,\ \eta}^\be$ is a constant
function and undefined if it is a one-to-one function.
\endproclaim
\1
\proclaim{(2.7) Remark}  It is immediate from Lemma 2.4 that for fixed
$i < \si$, and fixed $\eta < \tau < \k$, either
all the $F_i(\a,\ \tau)$ are defined or all the $F_i(\a,\ \tau)$ are
undefined, and similarly for the $F_i(\eta,\ \be)$.  Further, it is
immediate from Remark 2.5 that if all the $F_i(\a,\ \tau)$ are defined
then all the $F_i(\eta,\ \be)$ are defined iff the function
$F_i(\cdot,\ \tau)$ is constant (and, when both of these statements hold,
$F_i(\eta,\ \cdot)$ is also constant, with the same constant value),
and the analogous equivalence holds, starting from the hypothesis
that all the $F_i(\eta,\ \be)$ are defined.
\endproclaim
\1
\proclaim{(2.8) Definition}  For $\eta < \k$ and $\a \in B_\eta$,
we define $W_\a$ to be $\{F_i(\eta^\pr,\ \a)|i < \si,\ \eta^\pr < \eta\}
\cup \{F_i(\a,\ \tau)|i < \si,\ \eta < \tau > \k\}$.
\endproclaim
\1
Note that
for each $\eta < \k,\ \{ W_\a | \a \in B_\eta \}$ is a system of sets
of ordinals of power at most $\si + \k$.  We have stated in terms of
$\si + \k$
to emphasize that we are temporarily working without any
assumptions as to the order relationship between $\si$ and $\k$.
Thus, for all $\eta < \k$, we can find $B^*_\eta \subseteq B_\eta$,
with $card\ B^*_\eta = card\ B_\eta$ such that the
$(W_\a|\a \in B^*_\eta)$ form a $\Delta\text{-system}$ whose heart
we denote by $H_\eta$.  We also set $H = \bigcup\{H_\eta|\eta < \k\}$.
We further assume all of the following, for each $\eta < \k$:

\1\1
\roster
\item $(o.t.\ W_\a|\a \in B^*_\eta)$ has constant value, $o_\eta$;
for $\a \in B^*_\eta$, we let $(\ga^\a_\xi|\xi < o_\eta)$
be the increasing enumeration
of $W_\a$,
\item there is fixed $a_\eta \subseteq o_\eta$ such that
for all $\a \in B^*_\eta,\ a_\eta = \{\xi < o_\eta|\ga^\a_\xi \in H_\eta\}$,
\item there is fixed $b_\eta \subseteq (\si\times\eta\times o_\eta) \cup
(\si\times(\k \setminus (\eta + 1))\times o_\eta)$ such that for all $\a \in
B^*_\eta$ and all $(i,\ \nu,\ \xi),\ (i,\ \nu,\ \xi) \in b_\eta$ iff
either $(\nu < \eta$ and $\ga^\a_\xi = F_i(\nu,\ \a))$ or $(\eta < \nu$ and
$\ga^\a_\xi = F_i(\a,\ \nu))$.
\endroster
\1
\proclaim{(2.9) Lemma} If $\a_k \in B^*_{\eta_k},\
\be_k \in B^*_{\tau_k},\ \eta_k < \tau_k < \k,\
k = 0,\ 1$, and $\{ \a_0,\ \be_0\} \neq \{ \a_1,\ \be_1 \}$,
then $A(\a_0,\ \be_0) \cap A(\a_1,\ \be_1) \subseteq H$.
\endproclaim
\demo{Proof}  Suppose that $F_i(\a_0,\ \be_0) = F_j(\a_1,\ \be_1)$.
Let $\a^* \in B^*_{\eta_1},\ \a^* \not\in \{\a_0,\ \a_1\}$ and
let $\be^* \in B^*_{\tau_1},\ \be^* \not\in \{\be_0,\ \be_1\}$.
Then, by canonicity, $F_i(\a_0,\ \be_0) = F_j(\a^*,\ \be_1)$, so
$F_j(\a_1,\ \be_1) = F_j(\a^*,\ \be_1)$, which means that
$F_i(\a_0,\ \be_0) \in W_{\be_1}$.  By a similar argument,
$F_i(\a_0,\ \be_0) = F_j(\a_1,\ \be^*) = F_j(\a^*,\ \be^*) \in
W_{\be^*}$, and then, since $F_i(\a_0,\ \be_0) \in W_{\be_1} \cap W_{\be^*},\
F_i(\a_0,\ \be_0) \in H_{\tau_1} \subseteq H$, as required.
\enddemo
\1
Let $\bold P_0 = \bold P|(H \cup dom\ p)$
and let $V^{\pr} = V^{\bold P_0}$.
Note that all our hypotheses on $V$ still hold in $V^\pr$ and
$V^{\bold P} = (V^\pr)^{\bold Q}$, where, in $V^\pr, \bold Q \cong
\bold P$.  Thus, we can first force with $\bold P_0$ without changing
anything relevant; therefore, we can assume that
$H,\ p = \emptyset$, which
we do, from here on.  By Lemma 2.9, this, of course, guarantees that

\1
$$ (*)\ \ \text{For } \a_i,\ \be_i \text{ as in Lemma 2.9, } A(\a_0,\ \be_0)
\cap A(\a_1,\ \be_1) = \emptyset.$$
\1
Now choose $\a_\eta \in B^*_\eta$ for $\eta < \k$.
\1

It is at this point that the proof of Theorem 4*, in \S 3, will begin
to diverge.
Here, we will assume that $p$ also forces
that $\bold c$ has no homogeneous green set of order type
$\th + 1$ and we will show that $p$ forces that
$\bold c$ has a homogeneous red set of power $\la$, while in
\S 3, in the proof of Theorem 4*,
our treatment of the colors will be more \lq\lq symmetrical".
However, the remainder of the argument, here, will be similar
quite similar in spirit to the argument in Case 2 in \S 3, below.

Recall that here, we have already argued that, in
$V^{\bold P},\ \k \rightarrow (\k,\ \th + 1)^2$.
Thus, in $V^{\bold P}$,
there must be $S \subseteq \k$ of power $\k$ such that
$\{\a_\eta|\eta \in S\}$ is homogeneous red for
$\bold c|\{\a_\eta|\eta < \k\}$ (and therefore also for $\bold c).$

So, let $\bold S$ be a $\bold P \text{-name and } q \in P,\ p \leq q$ be
such that $q$ forces that $\{\a_\eta|\eta \in \bold S\}$ is
homogeneous red for $\bold c$ and that $\bold S$ has power $\k$.
Then, in $V$, there are $S \subseteq \k$, and for $\eta \in S,\
q_\eta \in P,\ q \leq q_\eta$ such that $q_\eta$ forces that
$\a_\eta \in \bold S$.

We may assume, without loss of generality, that the $(q_\eta|\eta \in S)$
form a $\Delta\text{-system}$ with heart $q$ (by which we mean that
the $q_\eta$ are pairwise isomorphic as well).
Thus, the $q_\eta$, for $\eta \in S$ are pairwise compatible and whenever
$\eta < \tau$ are both in $S$ and $q_\eta,\ q_\tau \leq r \in P,\ r$
forces that $\bold c (\a_\eta,\ \a_\tau) =$ red.

Let $A^* = \bigcup\{A(\a_\eta,\ \a_\tau)|\eta < \tau, \text{ both in } S\}$.
We may also assume that for all $\eta \in S,\ dom\ q_\eta \subseteq A^*$.
This is because if this fails, then, letting $\overline{q}_\eta =
q_\eta|A^*$, whenever $\eta < \tau$ are both in $S$ and
$\overline{q}_\eta,\ \overline{q}_\tau \leq \overline{r} \in P,\
\overline{r}$ forces that $\bold c (\a_\eta,\ \a_\tau) =$ red, because
$\bold c (\a_\eta,\ \a_\tau)$ is a
$\bold P |A(\a_\eta,\ \a_\tau)\text{-name}$ and $r$ forces that
$\bold c (\a_\eta,\ \a_\tau) =$ red, where $r = \overline{r} \cup
(q_\eta \setminus \overline{q}_\eta) \cup
(q_\tau \setminus \overline{q}_\tau)$, and this is all that is required
for the rest of the argument.

Further, we can clearly thin out $S$ to obtain a subset, $S^\pr$,
also of power $\k$, such that
for $\tau \in S^\pr$, letting $\tau^\pr = min\ S^\pr
\setminus (\tau + 1),\ dom\ q_\tau \setminus dom\ q \subseteq
\bigcup\{A(\a_{\eta_1},\ \a_{\eta_2})|\eta_1 <
\eta_2 < \tau^\pr\text{ both in } S\}$.

Finally, for $\tau \in S^\pr$ and $\a \in B^*_\tau$, we make a
copy $q^\tau_\a$ of $q_i$, above $q$.  We do this by moving only
coordinates in the $A(\a_\eta,\ \a_\tau)$ and the $A(\a_\tau,\ \a_\eta)$
which are in $dom\ q_\tau \setminus dom\ q$.  We move these coordinates
according to the order-isomorphisms between the $A(\a_\eta,\ \a_\tau)$ and
the $A(\a_\eta,\ \a)$, and the order-isomorphisms between the
$A(\a_\tau,\ \a_\eta)$ and the $A(\a,\ \a_\eta)$.  Clearly by
$(\ast)$, above
and by the previous paragraph, this is well-defined.
Also, by Lemma 2.9 and $(\ast)$, the $q^\tau_\a$ are pairwise
compatible.

Further, arguing as in the paragraph immediately preceding
Lemma 2.2, it is easy to see that whenever $\eta < \tau$ are
both in $S^\pr,\ \a \in B^*_\eta,\ \be \in B^*_\tau$ and
$q^\eta_\a \cup q^\tau_\be \leq r,\ r$ forces
that $\bold c (\a,\ \be) = \text{ red}$.
Finally, clearly,
whenever $q \leq r \in P,\ r$ is incompatible with fewer than
$\la$ many of the $q_\a^\tau$ for $\tau \in S^\pr$ and $\a \in B^*_\tau$.
But then, letting $\bold G$ be the canonical $\bold P\text{-name}$ for
the generic, let $\bold Y$ be the following $\bold P\text{-name}$:
\1
$$\ \{\a|\text{ there is } \tau \in S^\pr \text{ such that } \a \in B^*_\tau
\text{ and } q_\a^\tau \in \bold G\}$$

\noindent
But then $q$ forces that $\bold Y$ has power $\la$ and is homogeneous
red for $\bold c$.  This concludes the proof of Theorem 3.
\enddemo
\1\1
\subheading{\S 3.  ADDING COHEN SETS ABOVE THE COFINALITY}
\1
Recall that in the Introduction we have already argued that
Theorem 4 follows from Theorem 4*.  Here, we will prove Theorem
4*, whose statement we recall.
\1\1
\proclaim{Theorem 4*}  Suppose that in  $V,\ \om < \k < \si =
\si^{<\si} < \la,\ \k = cf\ \la \leq \nu, \la$ is a strong limit cardinal
and $\bold P$ is the partial ordering for adding $\nu$ Cohen subsets of
$\si$.  Suppose, further, that $\zeta < \k$ and that $\k \rightarrow
(\k,\ \zeta)^2$.  Then, in $V^{\bold P},\ \la \rightarrow
(\la,\ \zeta)^2$.
\endproclaim
\1

\demo{Proof}  We carry over from \S 2 all the material
up to and including the choice of the $\a_\eta \in B^*_\eta$, for
$\eta < \k$, and in particular, (2.1) - (2.9),
except that here, the analogue of (!) of \S 2 is:
\1

$$\align (!!)\ \  &\text{ in } V^{\bold P}, \text{ for each }
\eta < \k \text{ which is either } 0\\
&\text{ or a successor ordinal, }
\la_\eta \rightarrow (\la_\eta,\ \zeta)^2.\endalign$$
\1
The argument for this exactly follows that for (!) in \S 2.
Also, as noted in the Introduction, it follows from
the hypotheses of the Theorem, that in $V^{\bold P},\
\k \rightarrow (\k,\ \zeta)^2$.
Once again, (!!) enables us to assume, without loss of generality,
that $p$ forces that that each $\De_\eta$ is homogeneous red for
$\bold c$.

Note that it is an easy consequence of Lemma 2.9 and
our assumption that $H = \emptyset$ that if $\eta < \tau < \k,\
q \in P|H_{\a_\eta},\ r \in P|H_{\a_\tau}$ then $q$ and $r$ are
compatible.  Recall that by the paragraph immediately preceding
$(\ast)$, of \S 2,
we are assuming that $p = \emptyset$.  Let $s \in P$.
We now argue, using the $\si\text{-completeness}$ of $\bold P$ and
the fact that $\si > \k$, that:
\1\1
\proclaim{Lemma 3.1}  In $V$, there is $(p_\eta|\eta < \k)$ such
that for all $\eta < \k,\ s|H_{\a_\eta} \leq p_\eta,\
dom\ p_\eta \subseteq H_{\a_\eta}$, and such that
\1
$$\align (\ast\ast)\ \ &\text{if } \eta < \tau < \k,\
x \in \{\text{red, green}\}\text{ either }\\
& s \cup p_\eta \cup p_\tau\text{ forces }
\bold c(\a_\eta,\ \a_\tau) = x\text{ or}\\
& \text{there is } q \in P|A(\a_\eta,\ \a_\tau) \text{ such that}
\endalign$$

\roster
\item $s \cup p_\eta \cup p_\tau \leq q$,
\item $q$ forces $\bold c (\a_\eta,\ \a_\tau) \neq x$,
\item $dom\ q \setminus (dom\ p_\eta \cup dom\ p_\tau) \subseteq
A(\a_\eta,\ \a_\tau) \setminus (H_{\a_\eta} \cup H_{\a_\tau})$.
\endroster
\endproclaim
\demo{Proof}  Let $((\eta_\gamma,\ \tau_\gamma)|\gamma < \k)$ enumerate
all the pairs $(\eta,\ \tau)$ with $\eta < \tau < \k$.
For $\gamma \leq \k$, we define $(p^\gamma_\eta|\eta < \k)$ by
recursion on $\gamma$
so that $p^\gamma_\eta \in P|H_{\a_\eta}$, and for all $\eta < \k$ and
all $\gamma_1 < \gamma_2 \leq \k,\ p^{\gamma_1}_\eta \leq p^{\gamma_2}_\eta$.

For $\eta < \k$, let $p^0_\eta = s|H_{\a_\eta}$
and for nonzero limit ordinals,
$\de < \k$, and $\eta < \k$,
let $p^\de_\eta = \bigcup\{p^\gamma_\eta|\gamma < \de\}$.  So, suppose
that $\gamma = \xi + 1$.  If $\eta \not\in  \{\eta_\xi,\ \tau_\xi\}$
we take $p^\gamma_\eta = p^\xi_\eta$.
We construct $p^\gamma_{\eta_\xi},\ p^\gamma_{\tau_\xi}$.
Let $\eta = \eta_\xi,\
\tau = \tau_\xi,\ \a = \a_\eta,\ \a^\pr = \a_\tau$, and let
$p(0) = p^\xi_\eta,\ p^\pr(0) = p^\xi_\tau$.  Identify red with $0$
and green with $1$.  We will have $p^\gamma_\eta = p(2),\ p^\gamma_\tau =
p^\pr(2)$, where we define $p(i),\ p^\pr(i),\ i = 1,\ 2$ by
the following two-stage recursion.
If $k = 0,\ 1$ and $p(k),\ p^\pr(k)$ are defined, and if
$s \cup p(k) \cup p^\pr(k)$ forces $\bold c (\a,\ \a^\pr) = k$, then we
set $p(k + 1) = p(k),\ p^\pr(k + 1) = p^\pr(k)$.  Otherwise,
choose $q \in P|A(\a,\ \a^\pr)$ such that $s \cup p(k) \cup p^\pr(k) \leq
q$ and such that $q$ forces $\bold c(\a,\ \a^\pr) = 1 - k$.
Finally, let $p(k + 1) = q|H_\a,\ p^\pr(k + 1) = q|H_{\a^\pr}$.

Clearly then, by construction, for $\eta < \k$, taking
$p_\eta = p^\k_\eta,\ p_\eta$ is as required.  This completes
the proof of the Lemma.
\enddemo\enddemo
\1
\proclaim{Remarks}  Although we have developed it for both colors,
we only use the machinery of $(\ast\ast)$ of Lemma 3.1 with
$x = \text{ red}$.
Also, in $(\ast\ast)$, if $s \cup p_\eta \cup p_\tau$ does not force
that $\bold c(\a_\eta,\ \a_\tau) = \text{ red}$, we choose
$q_{\eta,\ \tau}$ to be some $q$ whose existence is guaranteed by
$(\ast\ast)$.
\endproclaim
\1
Now, still working in $V$, we define $d:[\k]^2 \rightarrow
\{\text{red, green}\}$ by $d(\eta,\ \tau) = \text{ red iff }
s \cup p_\eta \cup p_\tau \text{ forces }
\bold c(\a_\eta,\ \a_\tau) =
\text{ red.}$  Now, in $V,\ \k \rightarrow (\k,\ \th + 1)^2$, so
either (Case 1) there is $Y \in [\k]^\zeta$ which is
homogeneous green for $d$, or (Case 2) there is $Y \in
[\k]^\k$ which is homogeneous red for $d$.  We show that in
Case 1, $s$ has an extension which forces that
there is a set of order type $\zeta$ which is homogeneous
green for $\bold c$, while in Case 2, $s$, itself, forces that
there is a set of power $\la$ which is homogeneous red for
$\bold c$.  Clearly this suffices, since then the empty condition
forces that $\la \rightarrow (\la,\ \zeta)^2$.
We consider the cases separately.
\1\1
\subheading{Case 1:  The Green Case}
\1
Let $Y \in [\k]^\zeta$ be homogeneous green for $d$.  For
$\eta < \tau$ both in $Y$, note that $q_{\eta,\ \tau}$ is defined,
since $d(\a_\eta,\ \a_\tau) = green$.  Set
\1

$$r = \bigcup\{q_{\eta,\ \tau}|\eta <
\tau \text{ both in } Y\}$$

Once we have argued that $r$ is a function, it will be clear
that $r \in P,\ s \leq r$ (since for any $\eta < \tau$ which
are both in $Y,\ s \leq q_{\eta,\ \tau}\text{)}$
and further that $r$ forces that $\{\a_\eta|\eta \in Y\}$
is homogeneous green for $\bold c$, since, again, whenever
$\eta < \tau$ are both in $Y,\ q_{\eta,\ \tau}$ forces
that $\bold c (\a_\eta,\ \a_\tau) = \text{ green.}$  But,
once again, it follows from the conjunction of Lemma 2.9
and $(\ast)$ that $r$ is a function.  This completes the
proof in Case 1.
\1\1
\subheading{Case 2:  The Red Case}
\1
As we already noted there, the last part of the
argument in \S 2 is quite similar in spirit
to the argument we shall give for this case.
Let $Y \in [\k]^\k$ be homogeneous red for $d$.
As in \S 2,
for $\eta \in Y$ and
$\a \in B^*_\eta$, let $p^\eta_\a = T(p_\eta)$,
where $T$ is the order isomorphism between
$H_{\a_\eta}$ and $H_\a$.  Once again, the
$p^\eta_\a\ (\eta \in Y,\ \a \in B^*_\eta)$ are
pairwise compatible, by Lemma 2.9 and $(\ast)$,
and whenever $\eta < \tau$ are both in $Y,\ \a \in
B^*_\eta,\ \be \in B^*_\tau$ and
$s \cup p^\eta_\a \cup
p^\tau_\be \leq q,\ q$ forces that $\bold c (\a,\ \be) =
\text{ red,}$ by the fact that $d(\eta,\ \tau) = \text{ red}$
and by the argument of the paragraph immediately preceding
Lemma 2.2.  Also, once again, for all $s \leq q \in P,\ q$
is incompatible with at most $\si$ of the $p^\eta_\a$.

Now, let $\bold G$ again be the canonical $\bold P\text{-name}
$ of the generic, and for $\eta \in Y$, let
$\bold X_\eta$ be the $\bold P\text{-name }
\{\a \in B^*_\eta|p^\eta_\a \in \bold G\}$.  Then,
since $card\ B^*_\eta > \si,\ s$ forces that
$card\ \bold X_\eta = card\ B^*_\eta$.  We conclude
by noting that by the previous paragraph,
$s$ also forces that \lq\lq if $\eta < \tau$ are both
in $Y,\ \a \in \bold X_\eta,\ \be \in \bold X_\tau$ then
$\bold c (\a,\ \be) = \text{red."}$.  In other words,
\lq\lq as promised", $s$ forces that
$\bigcup\{\bold X_\eta|\eta \in Y\}$
is homogeneous red for $\bold c$ and has power $\la$.
This concludes the proof of Case 2, and therefore of
Theorem 4*.
\1\1

\enddocument

\end